\documentclass[11pt]{article}

\usepackage{xspace,amssymb,amsmath,amscd,amsthm,epsfig,array}
\makeatletter
\newfont{\footsc}{cmcsc10 at 8truept}
\newfont{\footbf}{cmbx10 at 8truept}
\newfont{\footrm}{cmr10 at 10truept}
\renewcommand{\ps@plain}{%
\renewcommand{\@oddfoot}{\footsc 
  \hfil\footrm\thepage}}
\makeatother
\theoremstyle{plain}
\newtheorem{theorem}{Theorem}

\newtheorem{conjecture}[theorem]{Conjecture}
\theoremstyle{definition}

\DeclareMathOperator{\area}{area}
\DeclareMathOperator{\dinv}{dinv}
\DeclareMathOperator{\adj}{adj}
\DeclareMathOperator{\sgn}{sgn}
\DeclareMathOperator{\spin}{spin}
\DeclareMathOperator{\bounce}{bounce}
\DeclareMathOperator{\id}{id}
\DeclareMathOperator{\maj}{maj}
\DeclareMathOperator{\inv}{inv}
\DeclareMathOperator{\diag}{diag}
\newcommand{\mymat}[3]{{}_{#3}[{#1}]_{#2}}
\newcommand{\scprod}[2]{\langle #1, #2 \rangle}

\newcommand{\ferrers}[1]{\mathcal{F}({#1})}
\newcommand{\nabq}{\nabla_{q=1}}
\newcommand{\Hmuq}[1]{\tilde{H}_{#1}|_{q=1}}

\newcommand{\N}{\mathbb{N}}

\newcommand{\Q}{\mathbb{Q}}
\newcommand{\lndp}[1]{LNDP_{#1}}
\newcommand{\ndp}[1]{NDP_{#1}}
\newcommand{\ldpn}{LDP_n}
\newcommand{\lndparb}[2]{\lndp{#2}^{#1}}
\newcommand{\lndpm}[1]{\lndparb{m}{#1}}
\newcommand{\sumsb}[1]{\sum_{\substack{#1}}}
\newcommand{\boldt}{\boldsymbol{T}}
\newcommand{\LLT}{\mathrm{LLT}}
\newcommand{\ssyt}{\mathrm{SSYT}} 
\newcommand{\srht}{\mathrm{SRHT}}
\newcommand{\myfig}[3]{\begin{figure}[htbp]
\begin{center}
{\scalebox{#1}{\includegraphics{#2}}}\caption{#3}\label{fig:#2}
\end{center}
\end{figure}}

\title{Nested Quantum Dyck Paths and $\nabla(s_{\lambda})$}

\author{
Nicholas A. Loehr \\
\small Department of Mathematics\\
\small Virginia Tech \\
\small Blacksburg, VA 24061-0123 \\
\small \texttt{nick@math.wm.edu}  \\
 \and
Gregory S. Warrington \\
\small Department of Mathematics\\
\small Wake Forest University \\
\small Winston-Salem, NC 27109 \\
\small \texttt{warrings@wfu.edu} 
}

\begin{document}
\maketitle

\begin{abstract}
We conjecture a combinatorial formula for the monomial expansion
of the image of any Schur function under the Bergeron-Garsia nabla operator.
The formula involves nested labelled Dyck paths weighted by area and a suitable
``diagonal inversion'' statistic. Our model includes as special cases
many previous conjectures connecting the nabla operator to quantum 
lattice paths.  The combinatorics of the inverse Kostka matrix
leads to an elementary proof of our proposed formula when $q=1$.
We also outline a possible approach for proving all the extant
nabla conjectures that reduces everything to the construction of 
sign-reversing involutions on explicit collections of 
signed, weighted objects.  
\end{abstract}

\section{Introduction}
\label{sec:intro} 

The \emph{nabla operator} introduced by Francois Bergeron
and Adriano Garsia~\cite{nabla1} plays a fundamental role in the 
theory of symmetric functions and Macdonald polynomials. To define
this operator, let us first introduce some notation. We let
$\Lambda$ denote the ring of symmetric functions in the variables
$x_1,x_2,\ldots$ with coefficients in the field $\Q(q,t)$.
The vector space $\Lambda$ has many well-known bases, all indexed
by integer partitions. We will use the following bases of 
$\Lambda$ in this paper: the monomial symmetric functions 
$m_{\mu}$; the homogeneous symmetric functions $h_{\mu}$;
the elementary symmetric functions $e_{\mu}$; the power-sum
symmetric functions $p_{\mu}$; the Schur functions $s_{\mu}$;
and the modified Macdonald polynomials $\tilde{H}_{\mu}$. More
details may be found in the encyclopedic reference~\cite{macbook}.

The nabla operator is the unique $\Q(q,t)$-linear map on $\Lambda$ 
such that $\nabla(\tilde{H}_{\mu})=q^{n(\mu')}t^{n(\mu)}\tilde{H}_{\mu}$
for all partitions $\mu$, where $n(\mu)=\sum_{i\geq 1} (i-1)\mu_i$
and $\mu'$ is the transpose of $\mu$. Thus, the modified 
Macdonald polynomials are the eigenfunctions of the nabla operator.
From the combinatorial point of view, the nabla operator is
important because it encodes a wealth of information about
$q,t$-analogues of combinatorial objects such as lattice paths,
parking functions, and labelled forests. The connection to
combinatorics arises by considering the matrix of the linear
operator $\nabla$ relative to various bases for $\Lambda$.
Given any two bases $(b_{\lambda})$ and $(c_{\mu})$ of $\Lambda$
and any linear operator $T$ on $\Lambda$,
we write $\mymat{T}{(b_{\lambda})}{(c_{\mu})}$ to denote the
unique matrix of scalars $a_{\mu,\lambda}\in\Q(q,t)$ such that
\[ T(b_{\lambda})=\sum_{\mu} a_{\mu,\lambda}c_{\mu} \]
for all partitions $\lambda$. In particular, if $T=\nabla$
and $(c_{\mu}')$ is the dual basis for $c_{\mu}$ relative to the Hall 
inner product on $\Lambda$, it follows that
\[ a_{\mu,\lambda}=\scprod{\nabla(b_{\lambda})}{c_{\mu}'}. \] 
We often restrict consideration to the subspace $\Lambda^n$ of symmetric 
functions of degree $n$, so that the matrix in question is a finite
square matrix with rows and columns indexed by the partitions of $n$.

By definition, $\mymat{\nabla}{(\tilde{H}_{\mu})}{(\tilde{H}_{\lambda})}$ is 
a diagonal matrix with diagonal entries $T_{\mu}=q^{n(\mu')}t^{n(\mu)}$.
For other choices of the input and output bases, one obtains other
$\Q(q,t)$-matrices representing the nabla operator. Remarkably, the
entries in these matrices are often \emph{polynomials} in $q$ and $t$
with integer coefficients all of like sign; i.e., we often have
$a_{\mu,\lambda}\in\pm\N[q,t]$. Whenever this occurs, one can seek
combinatorial interpretations for various entries $a_{\mu,\lambda}$
as sums of suitable signed, weighted objects. Such interpretations
have been sought after, conjectured, and (in some cases) proved
by many different authors. Table~\ref{tab:nabla} gives a list (not
necessarily exhaustive) of some recent research efforts in this area.

\begin{table}\setlength{\extrarowheight}{0.05cm}
\begin{center}
{\small
\begin{tabular}{|l|l|l|l|} \hline
Algebraic Object & Combinatorial Model & Conjectured by: & Proved by: 
\\\hline\hline
$\scprod{\nabla(e_n)}{s_{1^n}}$ & $q,t$-Dyck paths & 
Haglund~\cite{haglundbounce} (cf. ~\cite{remarkable}) 
& Garsia,Haglund~\cite{nablaproof1,nablaproof2}\\\hline
$\scprod{\nabla(e_n)}{h_{1^n}}$ & $q,t$-parking functions &
Haglund,Loehr~\cite{parkconj,mypark} & open \\\hline
$\scprod{\nabla(e_n)}{e_dh_{n-d}}$ & $q,t$-Schr\"oder paths &
Egge,Haglund, & Haglund~\cite{schrproof}  \\
& & Kremer,Killpatrick~\cite{qtschr} & \\\hline
$\scprod{\nabla^m(e_n)}{s_{1^n}}$ ($m>1$) & $m$-Dyck paths & 
Loehr~\cite{mconj} & open \\\hline
$\scprod{\nabla^m(e_n)}{h_{1^n}}$ ($m>1$) & labelled $m$-Dyck paths & 
Loehr,Remmel~\cite{mhconj} & open \\\hline 
$\scprod{\nabla^m(e_n)}{h_{\mu}}$ ($m\geq 1$) 
& labelled $m$-Dyck paths & Haglund,Haiman,Loehr, & open \\
 & & Remmel,Ulyanov~\cite{hhlru} & \\\hline 
$\scprod{\nabla(p_n)}{s_{1^n}}$ & $q,t$-square paths & 
Loehr,Warrington~\cite{sqconj} & Can,Loehr~\cite{sqproof} \\\hline
$\scprod{\nabla(p_n)}{h_{\mu}}$ & labelled square paths & 
Loehr,Warrington~\cite{sqconj} & open \\\hline 
$\scprod{\nabq(s_{\lambda/\nu})}{s_{\mu}}$ & digraphs 
& Lenart~\cite{lenart} & Lenart~\cite{lenart} \\\hline 
$\scprod{\tilde{H}_{\mu}}{h_{\mu}}$ & fillings of $\ferrers{\mu}$ &
Haglund~\cite{mac-conj} & Haglund,Haiman, \\
 & & & Loehr~\cite{mac-proof,mac-proof2}\\\hline
\end{tabular}
\caption{Summary of research on the combinatorics of the nabla operator.}
\label{tab:nabla}
}
\end{center}
\end{table}

Each of the conjectures mentioned in Table~\ref{tab:nabla} gives
only \emph{partial} information about the nabla operator. For example,
the Garsia-Haglund $q,t$-Catalan Theorem establishes a combinatorial
interpretation for just one of the coefficients in the matrix
$\mymat{\nabla}{(s_{\lambda})}{(s_{\mu})}$, namely
$a_{s_{1^n},s_{1^n}}=\scprod{\nabla(e_n)}{s_{1^n}}$.
The main conjecture in~\cite{hhlru} extends this result to a combinatorial
interpretation for the monomial expansion of $\nabla(e_n)$, but
this still only yields information about one column of the matrix
$\mymat{\nabla}{(s_{\lambda})}{(m_{\mu})}$.  Our goal in this paper 
is to present a new conjecture that gives a combinatorial interpretation
for \emph{every} entry in the matrix $\mymat{\nabla}{(s_{\lambda})}
{(m_{\mu})}$.  We shall see that this conjecture unifies and clarifies
the partial conjectures mentioned in Table~\ref{tab:nabla}.

In Section~\ref{sec:conjecture} of this paper, we describe our
conjectured combinatorial model for the monomial expansion of
$\nabla(s_{\lambda})$ and explain some connections to the
more specialized conjectures in Table~\ref{tab:nabla}.
In Section~\ref{sec:proof}, we give a proof of our conjecture
when $q=1$; the proof relies heavily on the combinatorics
of the inverse Kostka matrix $K^{-1}=\mymat{\id}{(m_{\lambda})}{(s_{\mu})}$.
In Section~\ref{sec:agenda}, we outline a combinatorial approach
that, if implemented, would prove the full conjecture. The proof method 
suggested in this final section, which relies heavily on the recently 
discovered combinatorial interpretation for modified Macdonald polynomials,
reduces all the extant nabla conjectures to the problem of defining
sign-reversing involutions on certain explicit collections of
signed, weighted objects.

\section{The Combinatorial Model} 
\label{sec:conjecture}

This section presents our conjectured formula for
the monomial expansion of $\nabla(s_{\lambda})$.
To prepare for this formula, we must first review
the known combinatorial interpretation for
$\scprod{\nabla(e_n)}{s_{1^n}}$ and the conjectured
interpretation for the monomial expansion of $\nabla(e_n)$.

\subsection{Quantum Dyck Paths}
\label{subsec:q-dyck}

Fix a positive integer $n$. A \emph{Dyck sequence} of length
$n$ is a list $g=(g_0,g_1,\ldots,g_{n-1})$ of nonnegative
integers such that $g_0=0$ and $g_{i+1}\leq g_i+1$ for all $i<n-1$.
The \emph{area} of a Dyck sequence is $\area(g)=\sum_{i=0}^{n-1} g_i$.
A \emph{diagonal inversion} of a Dyck sequence is a pair of
indices $i<j$ such that $(g_i-g_j)\in\{0,1\}$. We let
$\dinv(g)$ be the number of diagonal inversions of $g$.
Given any logical statement $P$, let $\chi(P)=1$ if $P$
is true, and $\chi(P)=0$ if $P$ is false. Then we can write
\[ \dinv(g)=\sum_{0\leq i<j<n} \chi(g_i-g_j\in\{0,1\}). \]
For example, $g=(0,0,1,2,0,1,1,2,3,1)$ is a Dyck sequence of
length $10$ with $\area(g)=11$ and $\dinv(g)=15$.

Dyck sequences correspond naturally to \emph{Dyck paths},
which are lattice paths from $(0,0)$ to $(n,n)$ consisting of $n$ unit 
north steps and $n$ unit east steps that never go below the line $y=x$.
We convert a Dyck sequence to a Dyck path by drawing $g_i$ complete
lattice squares to the left of the line $y=x$ in the $i$'th row
from the bottom, and following the north and west boundary of these
squares to obtain a lattice path. Then $\area(g)$ is the number of
squares between the path and the line $y=x$. The statistic $\dinv(g)$
counts pairs of cells lying immediately right of north steps in the
path, such that the cells are either on the same diagonal,
or such that the lower square lies one diagonal to the left of
the upper square. The diagonal inversion statistic was proposed
by Haiman in connection with the Garsia-Haiman $q,t$-Catalan
numbers. Haglund had previously defined another statistic on Dyck paths, called
the \emph{bounce score}~\cite{haglundbounce}.  Garsia and Haglund proved that
\[ \scprod{\nabla(e_n)}{s_{1^n}}=\sum_{g\in DS_n} t^{\area(g)}q^{\dinv(g)}
   =\sum_{\pi\in DP_n} t^{\bounce(\pi)}q^{\area(\pi)}, \]
where $DS_n$ is the set of Dyck sequences of order $n$,
and $DP_n$ is the set of Dyck paths of order $n$~\cite{nablaproof1,nablaproof2}.

In~\cite{hhlru}, the previous combinatorial formula was extended
to a conjectured formula for the monomial expansion of $\nabla(e_n)$.
To describe this extension, we consider pairs $(g,r)$, where
$g$ is a Dyck sequence of length $n$ and $r=(r_0,r_1,\ldots,r_{n-1})$
is a list of $n$ positive integers such that
$g_{i+1}=g_i+1$ implies $r_i<r_{i+1}$. 
Let $LDP_n$ be the set of all such pairs.  
We visualize an object $(g,r)\in LDP_n$ by drawing
the Dyck path $\pi$ associated to $g$ and labelling the $n$
north steps of $\pi$ (from bottom to top) with the labels
$r_0,r_1,\ldots,r_{n-1}$. The only restriction on the labels
is that the labels in each column must strictly increase reading upwards.
Now, define $\area(g,r)=\sum_{i=0}^{n-1} g_i=\area(g)$, and define
\[ \dinv(g,r)=\sum_{i<j} \chi(g_i-g_j=0\mbox{ and }r_i<r_j)
             +\sum_{i<j} \chi(g_i-g_j=1\mbox{ and }r_i>r_j). \]
Haglund, Haiman, Loehr, Remmel, and Ulyanov~\cite{hhlru} conjectured that
\begin{equation}\label{eq:hhlru}
 \nabla(e_n)=\sum_{(g,r)\in LDP_n}
     t^{\area(g,r)}q^{\dinv(g,r)}\prod_{i=0}^{n-1} x_{r_i}, 
\end{equation}
and they proved that this expression is symmetric in the $x_j$'s.
Because of this symmetry, the ``HHLRU conjecture'' is equivalent to the 
following assertion: for all $\mu=(\mu_1,\mu_2,\ldots)\vdash n$,
the $(\mu,(1^n))$-entry of $\mymat{\nabla}{(s_{\lambda})}{(m_{\mu})}$
is the sum of $t^{\area(g,r)}q^{\dinv(g,r)}$ over all
objects $(g,r)\in LDP_n$ such that the $r$-vector
contains $\mu_i$ copies of $i$ for all $i$.

\subsection{Combinatorial Model for $\scprod{\nabla(s_{\lambda})}{s_{1^n}}$}
\label{subsec:sign-conj}

We are going to conjecture a formula for the monomial expansion
of $\nabla(s_{\lambda})$. Before doing so, we describe a related
conjecture for the ``sign character'' $\scprod{\nabla(s_{\lambda})}{s_{1^n}}$,
where $\lambda$ is an arbitrary partition of $n$. The sign
character conjecture involves \emph{nested quantum Dyck paths},
while the full conjecture involves \emph{nested quantum labelled
Dyck paths}.

Our conjecture for the sign character has the form
\begin{equation}\label{eq:sign-conj}
 \scprod{\nabla(s_{\lambda})}{s_{1^n}}=
  \sgn(\lambda)\sum_{G\in NDP_{\lambda}} t^{\area(G)}q^{\dinv(G)}, 
\end{equation}
where $\sgn(\lambda)\in\{+1,-1\}$, $NDP_{\lambda}$ is a certain
collection of nested Dyck paths constructed from $\lambda$,
and $\area,\dinv:NDP_{\lambda}\rightarrow\N$ are suitable weight functions.
We will define the quantities $\sgn(\lambda)$, $NDP_{\lambda}$,
$\area$ and $\dinv$ in the context of a specific example.  

\begin{figure}
\begin{center}
\epsfig{file=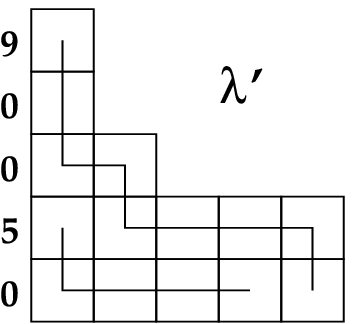}
\caption{Dissection of $\lambda'$ into border strips.}
\label{fig:border}
\end{center}
\end{figure}

Suppose $n=14$ and $\lambda=(5,3,2,2,2)$.
We begin by drawing the Ferrers diagram of the transposed partition
$\lambda'=(5,5,2,1,1)$. Next, we fill this diagram with ``rim hooks''
by repeatedly removing the entire northeast border of $\lambda'$,
as shown in Figure~\ref{fig:border}.  For $0\leq i<\lambda_1
=\ell(\lambda')$, let $n_i$ be the length of the hook that starts in 
the $i$'th row from the top of the diagram; let $n_i$ be zero
if there is no such hook. In our example, we have
\[ (n_0,n_1,n_2,n_3,n_4)=(9,0,0,5,0). \]
Define the \emph{spin} of $\lambda'$ to be the total number
of times a border hook crosses a horizontal boundary of a unit
square in the Ferrers diagram of $\lambda'$, and define
the sign $\sgn(\lambda)=(-1)^{\spin(\lambda')}$.
In our example, $\sgn(\lambda)=(-1)^{5}=-1$. 
We also define the \emph{dinv adjustment} by setting
\[ \adj(\lambda)=\sum_{i=0}^{\lambda_1-1} (\lambda_1-1-i)\chi(n_i>0). \]
This adjustment is the sum of the row indices in which 
the nonzero border hooks start, if we number the
rows $0,1,\ldots$ reading from bottom to top.
In our example, $\adj(\lambda)=1+4=5$.

\begin{figure}
\begin{center}
\epsfig{file=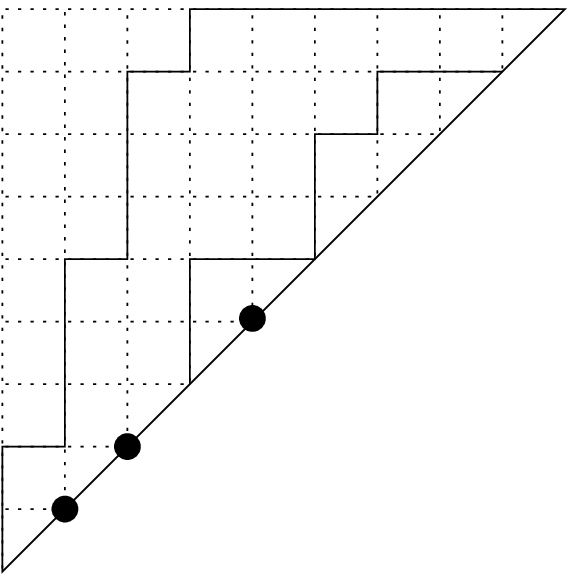}
\caption{Example of nested Dyck paths.}
\label{fig:nestpath}
\end{center}
\end{figure}

Next, we describe the collection of objects $NDP_{\lambda}$.
Let $l=\lambda_1=\ell(\lambda')$. We consider $l$-tuples
of lattice paths $\Pi=(\pi_0,\pi_1,\ldots,\pi_{l-1})$ such that $\pi_i$
is a lattice path from $(i,i)$ to $(i+n_i,i+n_i)$ consisting
of $n_i$ unit north steps and $n_i$ unit east steps that never
go strictly below the line $y=x$. If $n_i=0$ for some $i$,
then $\pi_i$ is a degenerate path consisting of a single vertex
at $(i,i)$. We say that $\Pi$ is \emph{nested} iff for all $i\neq j$,
no edge or vertex of $\pi_i$ coincides with any edge or vertex
of $\pi_j$. By definition, $NDP_{\lambda}$ consists of all such
$l$-tuples of nested Dyck paths. Note that degenerate paths are
important for determining nesting. Figure~\ref{fig:nestpath}
shows a typical element of $NDP_{(5,3,2,2,2)}$.

We can represent $\Pi=(\pi_0,\pi_1,\ldots,\pi_{l-1})\in NDP_{\lambda}$
by a ``Dyck configuration,'' which is the analogue of a Dyck sequence
in this setting. A Dyck configuration is an $l$-tuple of words
$G=(g^{(0)},g^{(1)},\ldots,g^{(l-1)})$, where $g^{(i)}$ is the Dyck sequence
for the Dyck path $\pi_i$. We choose the indexing of the letters
in these Dyck sequences to match the alignment of paths in the picture.
More precisely, for $i\leq a<i+n_i$, let $g^{(i)}_a$ be the number of complete 
lattice squares in the region bounded below by $y=a$, bounded above by 
$y=a+1$, bounded on the right by $y=x$, and bounded on the left by $\pi_i$; 
for all other values of $a$, $g^{(i)}_a$ is undefined. For example,
the Dyck configuration for the nested paths in Figure~\ref{fig:nestpath} is:
\[ G=\left(\begin{array}{cccccccccc}
g^{(0)}: & 0 & 1 & 1 & 2 & 3 & 3 & 4 & 5 & 5 \\
g^{(1)}: & . & . & . & . & . & . & . & . & . \\
g^{(2)}: & . & . & . & . & . & . & . & . & . \\
g^{(3)}: & . & . & . & 0 & 1 & 0 & 1 & 1 & . \\
g^{(4)}: & . & . & . & . & . & . & . & . & . \end{array}\right).  \]
Here, dots indicate positions where $g^{(i)}_a$ is undefined.

It is convenient to identify nested Dyck paths with the associated
Dyck configurations. Given a partition $\lambda=(\lambda_1,\lambda_2,
\ldots,\lambda_l)$ with associated border hooks of lengths
$n_1,n_2,\ldots,n_l$, an $l$-tuple $G=(g^{(i)}_a)$ belongs to $NDP_{\lambda}$ 
iff the following requirements are satisfied:
(i) for every $i$, $(g^{(i)}_a:i\leq a<i+n_i)$ is a Dyck sequence of 
length $n_i$; (ii) if $n_k=0$, then $g^{(i)}_k>0$ for all $i<k$
such that $g^{(i)}_k$ is defined; 
(iii) if $n_k>0$ and $i<k$, then $g^{(i)}_j>g^{(k)}_j$ for all $j$
such that both sides are defined; and (iv) if $n_k>0$ and $i<k$,
then $g^{(i)}_j>g^{(k)}_{j-1}+1$ for all $j$ such that both sides are defined.
Conditions (iii) and (iv) express the nesting requirement for two
nontrivial paths in terms of the Dyck sequences; condition (ii)
expresses the nesting requirement when the inner path has length zero.

Given $G\in NDP_{\lambda}$, the \emph{area} of $G$ is the sum of
the areas of the Dyck paths comprising $G$:
\[ \area(G)=\sum_{i=0}^{l-1}\sum_{i\leq j<i+n_i} g^{(i)}_j. \]
Note that lattice squares in the picture that appear inside
multiple Dyck paths are counted multiple times in the area statistic.
Next, the \emph{diagonal inversion statistic} for $G$ is
\begin{eqnarray*}
\dinv(G) &=& \adj(\lambda)+\sum_{a,b,u,v}
\chi(g^{(u)}_a-g^{(v)}_b=1)\chi(a\leq b) \\ & &
+\sum_{a,b,u,v} \chi(g^{(u)}_a-g^{(v)}_b=0)
\chi((a<b)\mbox{ or }(a=b\mbox{ and }u<v)).
\end{eqnarray*}
In these sums, we consider all possible choices
of $a,b,u,v$ such that $g^{(u)}_a$ and $g^{(v)}_b$
are both defined.  For the example shown in Figure~\ref{fig:nestpath},
we have $\area(G)=27$ and $\dinv(G)=28$.

Our conjectured formula 
\[ \scprod{\nabla(s_{\lambda})}{s_{1^n}}=
  \sgn(\lambda)\sum_{G\in NDP_{\lambda}} t^{\area(G)}q^{\dinv(G)} \]
has been verified by computer for all partitions $\lambda$
of all integers $n\leq 9$.

\subsection{Combinatorial Model for 
            the Monomial Expansion of $\nabla(s_{\lambda})$}
\label{subsec:combmodel}

Roughly speaking, our main conjecture for $\nabla(s_{\lambda})$ is
obtained by ``adding labels'' to the sign character formula introduced
in the previous subsection.  This is done in the spirit of the shuffle
conjecture for $\nabla(e_n)$ of~\cite{hhlru}.  For any partition
$\lambda$, we associate a collection of labelled nested Dyck paths
which we abbreviate $\lndp{\lambda}$.  Note that ``nested'' in this
case will be a slightly weaker notion than that used in the definition
of $\ndp{\lambda}$.

Let $l = \lambda_1 = \ell(\lambda')$.  An element of $\lndp{\lambda}$
consists of a pair $(G,R)$.  Here, $G$ is an $l$-tuple
$(g^{(0)},g^{(1)},\ldots,g^{(l-1)})$ where each $g^{(j)}$
encodes the Dyck sequence of some path $\pi_j$ of length $n_j$ among
the entries $(g_j^{(j)},g_{j+1}^{(j)},\ldots,g_{j+n_j-1}^{(j)})$.
The $l$-tuple $R = (r^{(0)},r^{(1)},\ldots,r^{(l-1)})$ is a list of
labels.  For all $j$, the length of $r^{(j)}$ equals the length of
$g^{(j)}$ (i.e., $r^{(j)} =
(r_j^{(j)},r_{j+1}^{(j)},\ldots,r_{j+n_j-1}^{(j)})$) and $r^{(j)}\in
(\mathbb{Z}_{>0})^{n_j}$.  Together, $G$ and $R$ are subject to the
following conditions: 
\begin{enumerate}
\item If $g_{i+1}^{(j)} = g_{i}^{(j)} + 1$, then
  $r^{(j)}_{i} < r^{(j)}_{i+1}$.\label{lndp:cond2}
\item The value $g_i^{(j)}$ is undefined or greater than zero for all
  $j < i \leq l-1$.\label{lndp:cond3}
\item For all $a$ and all $j < k$, either one of $g_a^{(j)}$ or
  $g_{a-1}^{(k)}$ is undefined, or $g_a^{(j)} > g_{a-1}^{(k)}$.\label{lndp:cond4}
\item For all $a$ and all $j<k$, 
if $g_a^{(j)}$ and $g_{a-1}^{(k)}$ are defined with $g_a^{(j)} =
  g_{a-1}^{(k)}+1$, then $r_a^{(j)} \leq r_{a-1}^{(k)}$.\label{lndp:cond5}
\end{enumerate}
The first condition states that every path $\pi_j$ is matched up with
a label vector $r^{(j)}$ that strictly increases up the columns of the
path.  The remaining conditions imply that no path encounters the
start of any other (even zero-length) path; that the paths are weakly
nested with no shared east steps; and that for a given
column, no larger label in the row directly above belongs to a
lower-indexed path.

Given $(G,R)\in\lndp{\lambda}$, the \emph{area} of $(G,R)$ is (as
before) the sum of the areas of the Dyck paths comprising G:
\begin{equation*}
  \area(G) = \sum_{i=0}^{l-1} \sum_{i\leq j < i + n_i} g^{(i)}_j.
\end{equation*}
The \emph{diagonal inversion statistic} for $(G,R)$ simply
incorporates the labels into the two summations:
\begin{equation}\label{eq:Gdinv}
\begin{aligned}
  \dinv(G,R) = \adj(\lambda) &+
  \sumsb{u,v,a,b}\chi(g^{(u)}_a-g^{(v)}_b=1)\chi(r^{(u)}_a > r^{(v)}_b)\chi(a\leq b) \\
  &+ \sumsb{u,v,a,b} \chi(g^{(u)}_a-g^{(v)}_b=0)
  \chi(r^{(u)}_a < r^{(v)}_b)\chi(a<b\text{ or }(a=b\text{ and }u<v)).
\end{aligned}
\end{equation}
In these sums, we consider all possible choices of $a,b,u,v$ such that
$g^{(u)}_a$ and $g^{(v)}_b$ are both defined.  
Set $x_R=\prod_{u,a} x_{r^{(u)}_a}$.
For the example shown in Figure~\ref{fig:lndpex},
we compute the coefficient of $m_{2,1,1}$ in $\nabla(s_{2,2})$ as
$-t^2q^2(2+t+q)$.  (Note that in this example, $\adj(\lambda) = 1$.)

\myfig{.5}{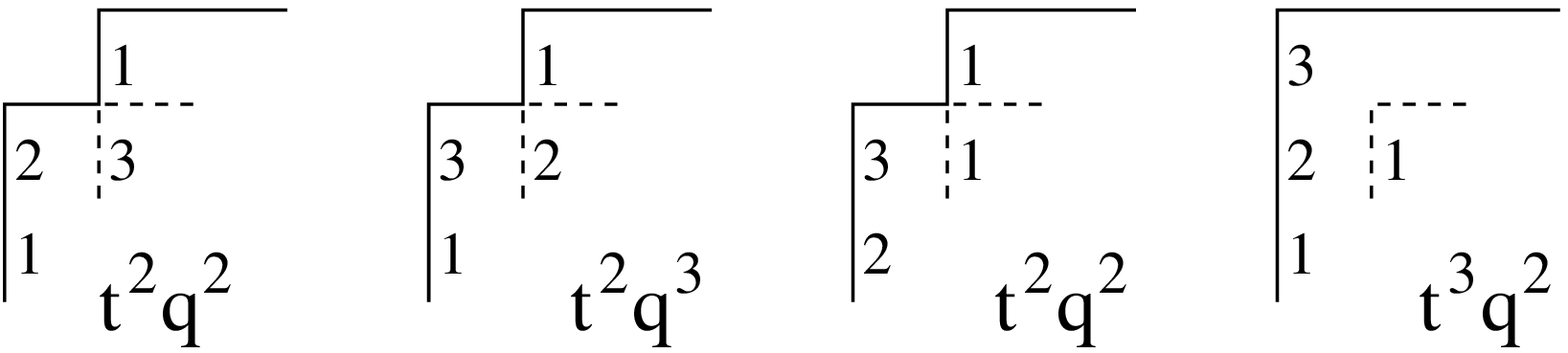}{Computation of the coefficient of $m_{2,1,1}$ in $\nabla(s_{2,2})$.}

In Figure~\ref{fig:lndptrip} we illustrate the pair $(G,R)$ where 
\begin{equation}\label{eq:lndpG}
   G=
  \begin{pmatrix}
    g^{(0)}: & 0 & 1 & 2 & 3 & 3 & 3 & 4 & 3 \\
    g^{(1)}: & . & 0 & 1 & 1 & 2 & 3 & 2 & . \\
    g^{(2)}: & . & . & 0 & 1 & 2 & . & . & . \\
  \end{pmatrix}
\end{equation}
and
\begin{equation*}
 R= \begin{pmatrix}
    r^{(0)}: & 2 & 6 &10 &15 &11 & 3 & 5 & 4 \\
    r^{(1)}: & . & 1 & 17 & 7 & 8 &14 &12 & . \\
    r^{(2)}: & . & . & 9 & 13 & 16 & . & . & . \\
  \end{pmatrix}.
\end{equation*}
We have $\area(G,R) = 31$ and $\dinv(G,R) = 24$.  

\myfig{.7}{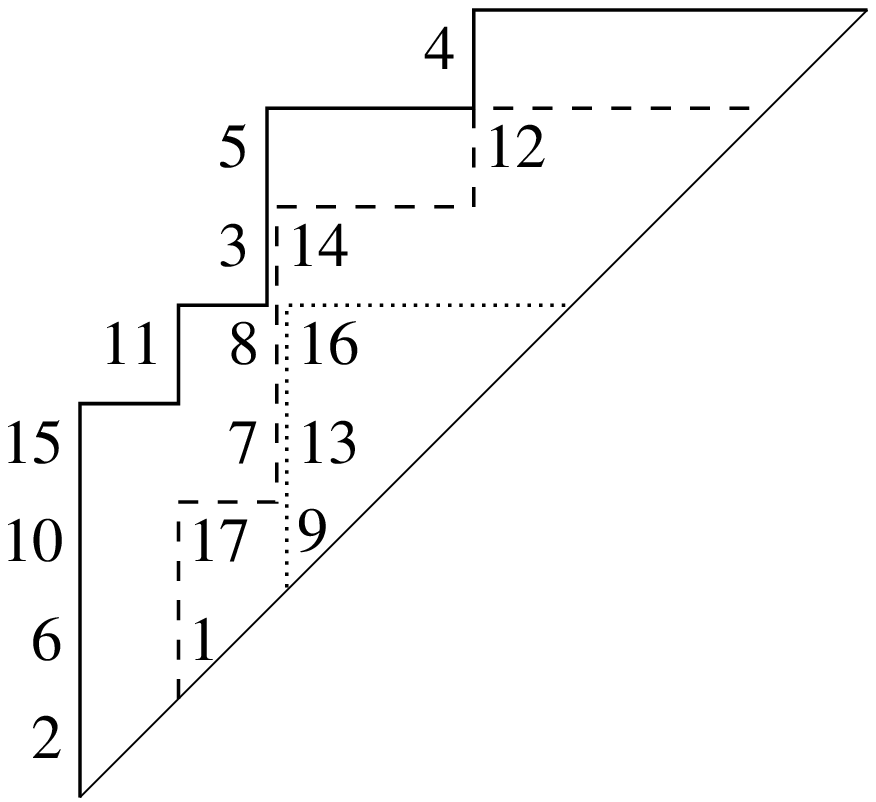}{One term contributing to the coefficient of
  $s_{1^{17}} = m_{1^{17}}$ in $\nabla(s_{3^5,2})$.}

\begin{conjecture}\label{conj:Main}
  For any partition $\lambda$, 
  \begin{equation}\label{eq:full-conj}
    \nabla(s_{\lambda})= \sgn(\lambda)
    \sum_{(G,R)\in \lndp{\lambda}} t^{\area(G,R)}q^{\dinv(G,R)}x_{R}.
  \end{equation}
\end{conjecture}

\begin{theorem}
  Our conjectured formula for $\nabla(s_\lambda)$ (i.e., the right
  side of \eqref{eq:full-conj}), is a symmetric function of the $x_j$'s.
\end{theorem}
\begin{proof}
  We prove this by expressing our summation in question as a weighted
  linear combination of the \emph{Lascoux-Leclerc-Thibon (LLT)
    polynomials} (introduced in~\cite{LLT}), which are known to be
  symmetric.

  The LLT polynomials can be defined combinatorially as follows.  Let
  $\Gamma$ be a $k$-tuple of skew shapes
  $(\lambda^1/\nu^1,\lambda^2/\nu^2,\ldots,\lambda^k/\nu^k)$ with $n$
  total boxes.  For each such $\Gamma$, we write $\ssyt_{\Gamma}^N$ to
  denote the set of $k$-tuples of semistandard Young tableaux $\boldt =
  (T_1,T_2,\ldots,T_k)$ such that for each $i$,
  \begin{enumerate}
  \item $T_i$ has shape $\lambda^i/\nu^i$.
  \item $T_i$ has entries from $\{x_1,x_2,\ldots,x_N\}$.
  \end{enumerate}
  We denote the content of $\boldt$ by $x_{\boldt} = 
      \prod_i\prod_{c\in\lambda^i/\nu^i} x_{T_i(c)}$.
  Given some $\boldt\in \ssyt^N_\Gamma$, we define a \emph{diagonal
    inversion} statistic as follows.  For a cell $c = (i,j)$, we
  define the \emph{diagonal} of $c$ to be $\diag(c) = j - i$.
  Then set
  \begin{equation}\label{eq:lltdinv}
    \begin{aligned}
    \dinv(\boldt) = \sumsb{k < l\\c\in T_k,\ d\in T_l}
    [&\chi(\diag(d) - \diag(c) = 1 \text{ and } T(c) < T(d))\\
     &+ \chi(\diag(c) - \diag(d) = 0 \text{ and } T(c) > T(d))].
    \end{aligned}
  \end{equation}
  The LLT polynomials are defined as
  \begin{equation*}
    \LLT^N_\Gamma(x_1,x_2,\ldots,x_N) = 
    \sum_{\boldt\in\ssyt^N_\Gamma} q^{\dinv(\boldt)} x_{\boldt}.
  \end{equation*}
  It is proved in~\cite{mac-proof,LLT} that each
  $\LLT^N_{\Gamma}(x_1,x_2,\ldots,x_N)$ is a symmetric polynomial in
  the $x_i$'s.  By taking inverse limits in the usual manner, we
  obtain symmetric polynomials $\LLT_{\Gamma} \in \Lambda^n$.

  We now explore the relationship between elements of $\lndp{\lambda}$
  and the LLT polynomials.  We describe the correspondence
  via an example.  Consider the pair $(G,R)$ illustrated in
  Figure~\ref{fig:lndptrip}.  Reading up the columns from right to
  left, we encounter five contiguous multisets of north steps of
  varying lengths.  The configuration $G$ will thereby be associated
  with a tuple of skew shapes $\Gamma(G) =
  (\lambda^1/\nu^1,\lambda^2/\nu^2,\ldots,\lambda^5/\nu^5)$ of sizes
  $2$, $8$, $2$, $1$ and $4$, respectively.  The third through fifth
  groups, each consisting of north steps from a single path, will
  yield skew shapes that are columns of the appropriate heights.  In
  particular, we get the shapes $(1^2)$, $(1^4)/(1^3)$ and $(1^4)$,
  respectively.  We have augmented as necessary each $\lambda^i$ and
  $\nu^i$ in equal amounts to ensure that any north step in $G$ going
  north from the line $y = x + b$ maps to a cell in $\Gamma(G)$ with
  diagonal equal to $b$.  When we have more than one path contributing
  to a given multiset of north steps, we proceed in an analogous
  manner.  However, in this case, the cells arising from the $j$-th
  path from the left are additionally shifted to the right and up by
  $j-1$ units.  (Notice that this is a diagonal-preserving shift.)
  Doing so yields $\lambda^1/\nu^1 = (2^4)/(2^3)$ and $\lambda^2/\nu^2
  = (3^5)/(3^2,1)$.
  
  The labels in $R$ accompany their respective north steps to give us
  the element $\boldt\in\ssyt^{17}_{\Gamma}$ illustrated in
  Figure~\ref{fig:tupleex}.  Note that we have aligned the skew shapes
  along the diagonals to facilitate computation of $\dinv$.

  \myfig{1}{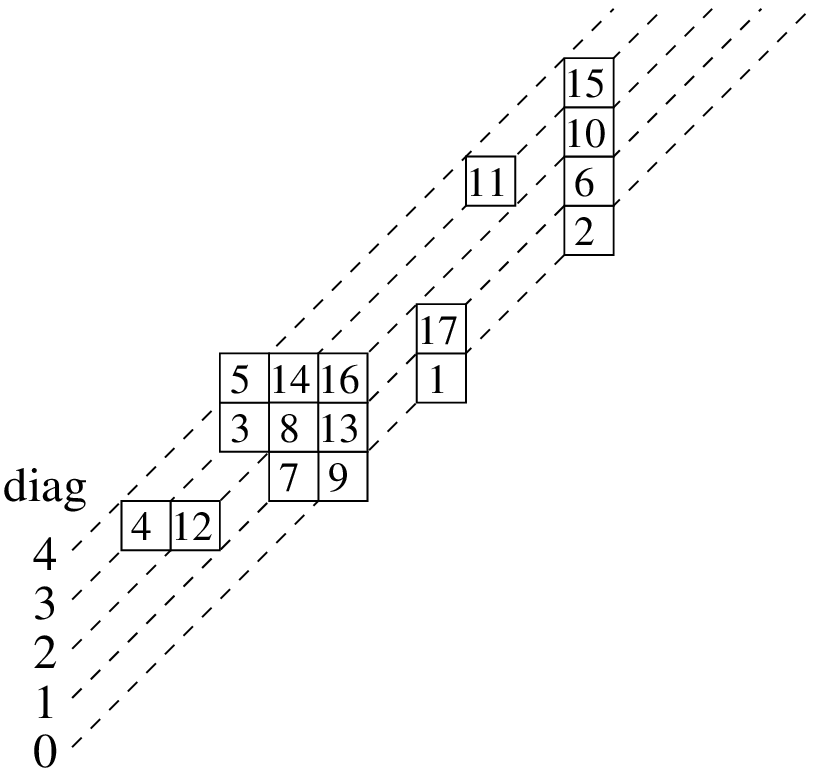}{$\boldt$ corresponding to $(G,R)$ of
    Figure~\ref{fig:lndptrip}.}

  The above correspondence gives a well-defined correspondence from a
  Dyck configuration $G$ to a tuple of shapes $\Gamma(G)$.  
  The correspondence can be modified to become
  invertible by the following two adjustments (which we do not make).  
  First, send all north steps in a given column to a (possibly 
  disconnected) skew shape.  Second, include empty skew shapes for 
  columns without any north steps. In any case, 
  by including labels, the correspondence between Dyck configurations
  and tuples of shapes extends to a map from elements $(G,R)\in
  \lndp{\lambda}$ to tuples of filled shapes $\boldt(G,R)$.  
  As is discussed below,
  these filled shapes are precisely the tuples $\boldt$ appearing in
  $\LLT_{\Gamma(G)}$.

  Our first claim is that the shapes are, in fact, skew shapes.  This
  follows from Conditions~\ref{lndp:cond3} and~\ref{lndp:cond4} in the
  definition of $\lndp{\lambda}$.  Our second claim is that the
  fillings of the skew shapes are semistandard.  That the entries
  increase up columns is the content of Condition~\ref{lndp:cond2}
  while rows are forced to weakly increase by Condition~\ref{lndp:cond5}.
  Furthermore, it is easily checked that any element
  $\boldt$ contributing to $\LLT_{\Gamma(G)}$ appears as the 
  image of $(G,R)$ for some labelling $R$.  It follows that 
  \begin{equation}\label{eq:gfcomp}
    \sum_{R:\ (G,R)\in\lndp{\lambda}} x_R = \sum_{\boldt\in \ssyt_{\Gamma(G)}} 
       x_{\boldt}.  
  \end{equation}

  We now consider the relationship between $\dinv(G,R)$ and
  $\dinv(\boldt(G,R))$.  Fix $(G,R)\in\lndp{\lambda}$
  corresponding to a $\boldt(G,R)\in\ssyt_{\Gamma(G)}$.  Suppose
  we have $u,v,a$ and $b$ such that $g_a^{(u)} - g_b^{(v)} =
  1$ and $r_a^{(u)} > r_b^{(v)}$.  Under our correspondence,
  $g_a^{(u)}$ and $g_b^{(v)}$ will correspond to cells $c$ and $d$
  (labelled $r_a^{(u)}$ and $r_b^{(v)}$, respectively) in $\boldt(G,R)
  = (T_1,T_2,\ldots,T_k)$.  We may write $c\in T_i$ and $d\in T_j$ for
  some $i,j$.  Since $a\leq b$ but $g_a^{(u)} > g_b^{(v)}$, we must
  have $i > j$.  Since diagonals are preserved, we see that $c$ and
  $d$ will contribute to the first summand in~\eqref{eq:lltdinv}.  A
  similar analysis of terms with $a < b$ in the second summation
  of~\eqref{eq:Gdinv} will yield the second summand of~\eqref{eq:lltdinv}.

  We are left to consider those quadruples $a,b,u,v$ for which $a = b$
  and $u < v$.  Such a quadruple will contribute to~\eqref{eq:Gdinv}
  exactly when $r_a^{(u)} < r_b^{(v)}$.  However, the corresponding
  pair of cells in $\boldt(G,R)$ will never contribute to
  $\dinv(\boldt(G,R))$ because they will lie in the same skew shape.
  Fortunately, for $(G,R)\in\lndp{\lambda}$ and such a quadruple, we
  will always have $r_a^{(u)} < r_b^{(v)}$: To see this, note that by
  Conditions~\ref{lndp:cond3} and~\ref{lndp:cond4}, $g_{a-1}^{(v)}$ is
  defined.  By Condition~\ref{lndp:cond5}, $r_a^{(u)} \leq
  r_{a-1}^{(v)}$.  Finally, $r_{a-1}^{(v)} < r_{a}^{(v)}$ by
  Condition~\ref{lndp:cond2}.  For brevity, write 
  \begin{equation*}
    \delta(\Gamma(G)) = \sum_{i=1}^k \delta(\lambda^i/\nu^i)
  \end{equation*}
  where we define $\delta$ for a skew shape by
  \begin{equation*}
    \delta(\lambda/\nu) = \sum_{(i,j)\in \lambda/\nu} 
	\sum_{a > 0} \chi((i+a,j+a)\in \lambda/\nu)
  \end{equation*} 
  Then
  \begin{equation}\label{eq:dinvcomp}
    \dinv(G,R) = \adj(\lambda) + \dinv(\boldt(G,R)) + n(\Gamma(G)).
  \end{equation}
  Combining~\eqref{eq:gfcomp} and~\eqref{eq:dinvcomp}, we conclude
  \begin{equation}\label{eq:corres}
    \sum_{R:\ (G,R)\in\lndp{\lambda}} q^{\dinv(G,R)} x_R = 
    q^{\adj(\lambda) + n(\Gamma(G))}\sum_{\boldt\in \ssyt_{\Gamma(G)}} 
    q^{\dinv(\boldt(G,R))} x_{\boldt}.
  \end{equation}
  Of course, the summation on the right side of~\eqref{eq:corres} is
  $\LLT_{\Gamma(G)}$.  So, if we allow $G$ to vary as well and sum
  over all such $G$, we get
  \begin{equation}\label{eq:fincorres}
    \sum_{(G,R)\in\lndp{\lambda}} t^{\area(G,R)} q^{\dinv(G,R)} x_R = 
    \sum_{G} t^{\area(G)} q^{\adj(\lambda) + n(\Gamma(G))} \LLT_{\Gamma(G)}.
  \end{equation}
  Since, as mentioned, each $\LLT_{\Gamma(G)}$ is a symmetric function
  in the $x_i$'s, it follows that the left side of
  \eqref{eq:fincorres} is as well.
\end{proof}

\subsection{Hook Shapes} 
\label{subsec:hook}

Let $\lambda = (a,1^{n-a})$ for some $a \geq 1$; so $\lambda' =
(n-a+1,1^{a-1})$.  When we fill the Ferrers diagram of $\lambda'$ with
rim hooks as described in Section~\ref{subsec:sign-conj}, we find that
there is only one nonzero rim hook. Thus, $(n_0,n_1,\ldots,n_{a-1})
=(n,0,\ldots,0)$, so $\sgn(\lambda) = (-1)^{a-1}$ and $\adj(\lambda) = a-1$.  
Furthermore, objects in $\lndp{\lambda}$ can be identified with those 
elements $(g,r)$ of $\ldpn$ for which $g_i > 0$ for $0 < i < a$.  
So when $\lambda$ is a hook, we obtain the following simplifications
of the main conjecture.

\begin{conjecture}\label{conj:hook}
  \begin{eqnarray*}
\scprod{\nabla(s_{(a,1^{n-a})})}{s_{1^n}} &=& (-q)^{a-1}
 \sum_{\substack{ g\in DP_n: \\ g_i>0\text{ for }0<i<a}} 
           t^{\area(g)}q^{\dinv(g)} \\
    \nabla(s_{a,1^{n-a}}) &=& (-q)^{a-1}\sum_{
\substack{(g,r)\in\ldpn: \\ g_i>0\text{ for }0<i<a}} 
t^{\area(g,r)}q^{\dinv(g,r)}x_r.
  \end{eqnarray*}
\end{conjecture}

\subsection{Trapezoidal Paths}
\label{subsec:trapz}
Conjecture~\ref{conj:Main} neatly explains some of the formulas
conjectured in~\cite{thesis,trap} regarding lattice paths in trapezoids.  
Consider the trapezoid bounded by the vertices $(0,k)$, $(0,k+n)$,
$(k,k)$ and $(k+n,k+n)$ for some $k,n$ with $k \geq 0$ and $n \geq 1$.
A \emph{trapezoidal lattice path} of type $(n,k)$ is a path from
$(0,k)$ to $(k+n,k+n)$ consisting of $n$ north steps and $k+n$ east
steps (all of length one) such that no vertex of the path lies
strictly below the line $y = x$.  The case of $k = 0$ is that of
Dyck paths.  Write $\mathcal{T}_{n,k}$ for the set of trapezoidal
lattice paths of type $(n,k)$.

Given a path $P\in \mathcal{T}_{n,k}$, define the
sequence $g(P) = (g_0,g_1,\ldots,g_{n-1})$ by taking $g_i$ to be the
number of unit squares in the strip bounded below by $y = k+i$, above
by $y=k+i+1$, on the left by $P$, and on the right by $y = x$.  Define 
$\area(P) = \sum_{i=0}^{n-1} g_i$  and 
\begin{equation*}
 \dinv(P) = \sum_{i < j} 
\chi(g_i-g_j\in \{0,1\}) + \sum_{i=0}^{n-1} \max(k-g_i,0).
\end{equation*}
It is conjectured in~\cite{thesis,trap} that 
$\sum_{P\in\mathcal{T}_{n,k}} t^{\area(P)}q^{\dinv(P)}$ 
is symmetric in $q$ and $t$.

We now show that this conjecture follows from
Conjecture~\ref{conj:Main}.  To this end, for $n,k \geq 0$ and $n + k
> 0$, define
\begin{equation*}
  \lambda(n,k) = \left(\left\lceil\frac{k+1}{2}\right
    \rceil^{\left\lfloor \frac{k+1}{2}\right\rfloor},1^n\right).
\end{equation*}
The partition $\lambda(n,k)$ has been defined so that when the special
rim hooks are placed in $\lambda(n,k)'$, the hooks are of lengths $n_0
= k+n$ and $n_i = k-2i$ for $1\leq i\leq \lfloor k/2\rfloor$.  Since the
successive length differences are only $2$ and the shortest path is of
length $0$ or $1$, this forces the path of length $n_i$ for $i > 0$ to
consist of $n_i$ north steps followed by $n_i$ east steps.  The
outermost path, of length $k+n$ has the single restriction that it must
begin with $k$ north steps.  It follows that the elements
$G\in\ndp{\lambda(n,k)}$ are in natural bijection with the $g(P)$ for
$P\in\mathcal{T}_{n,k}$ by sending $G =
(g^{(0)},g^{(1)},\ldots,g^{(\lfloor k/2\rfloor)})$ to
$g' = (g^{(0)}_k,g^{(0)}_{k+1},\ldots,g^{(0)}_{k+n-1})$.

To complete the proof, we need only examine how $\area(G)$
compares to $\area(g')$ and how $\dinv(G)$ compares to
$\dinv(g')$.  It is a simple computation to show that
\begin{equation*}
  \area(G) = \area(g') + \sum_{i=0}^{\lfloor k/2\rfloor} \binom{k-2i}{2}.
\end{equation*}
As for the diagonal inversion statistic, first note that by throwing
away all but the outermost path, we have lost $\adj(\lambda(n,k)) =
\binom{\lfloor k/2\rfloor + 1}{2}$.  We have also lost the
contributions to $\dinv(G)$ arising from interactions between any two
of the paths.  However, we gain $\sum_{i=k}^{k+n-1}
\max(k-g^{(0)}_i,0)$.  We leave it to the reader to show, in fact,
that these collectively give the appropriate difference.

\myfig{.6}{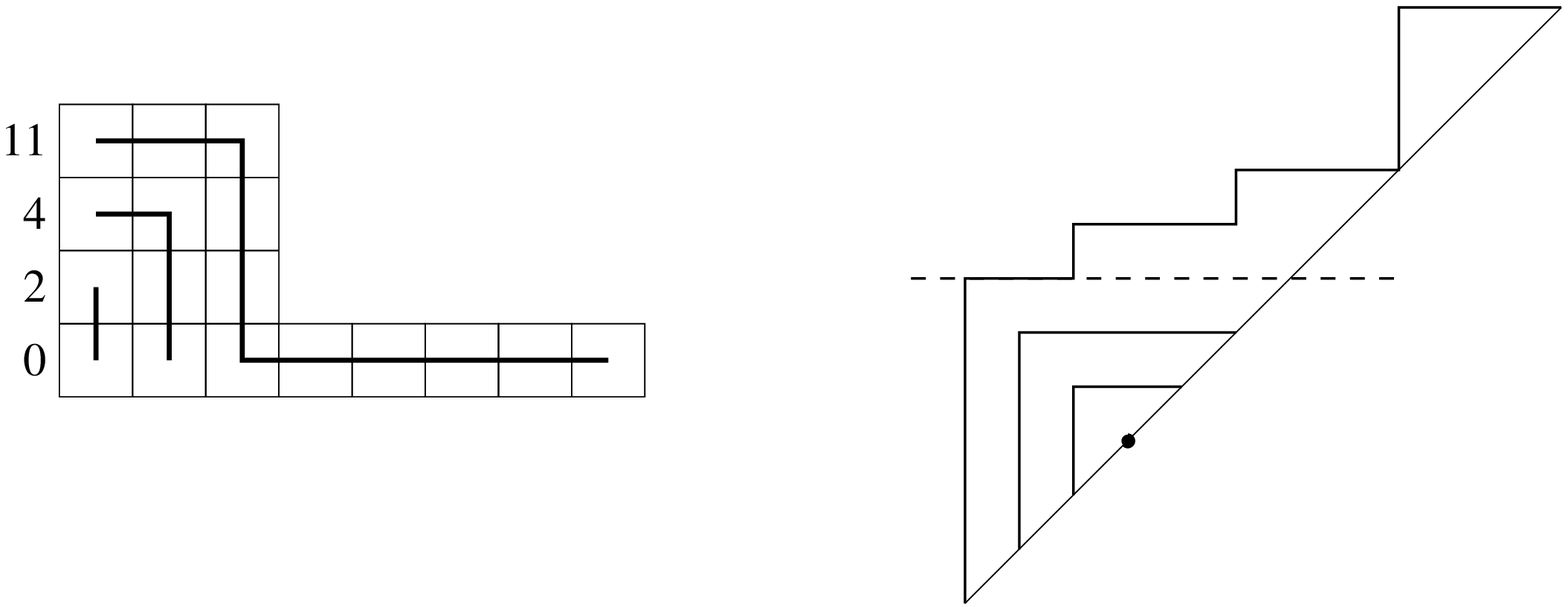}{Example of correspondence with trapezoidal paths.}

We illustrate an example in Figure~\ref{fig:trapii} for $k = 6$ and
$n=5$.  On the left is shown $\lambda(5,6)'$ along with its special
rim hooks of lengths $11$, $4$, $2$ and $0$.  On the right we show a
typical element of $\ndp{\lambda(5,6)}$.  Notice that by considering
the portion of the outermost path weakly above the dotted line, we
have an element of $\mathcal{T}_{5,6}$.  For this example,

\begin{equation*}
   G=
   \left(\begin{array}{cccccccccccc}
     g^{(0)}: & 0 & 1 & 2 & 3 & 4 & 5 & 4 & 2 & 0 & 1 & 2 \\
     g^{(1)}: & . & 0 & 1 & 2 & 3 & . & . & . & . & . & . \\
     g^{(2)}: & . & . & 0 & 1 & . & . & . & . & . & . & . \\
     g^{(3)}: & . & . & . & . & . & . & . & . & . & . & . \end{array}\right).
\end{equation*}
So $\area(G) = 31$, $\dinv(G) = 45$ and $\area(g') = 9$.  We now
check that $\dinv(G)-\dinv(g')$ is indeed $22 = 31-9$.  In going from
$G$ to $g'$, we lose the adjustment $\adj(\lambda(5,6)) = 6$.  We also
lose all of the interactions between $g^{(0)}$ and $g^{(i)}$ for
$i=1,2,3$ that contribute to $\dinv(G)$.  These account for a loss of
$37$ more.  However, we gain the sum $\sum_{i=0}^5 \max(6-g'_i,0) =
21$.  This gives us a difference of $22$ as desired.  

\subsection{Higher Powers of Nabla}
\label{subsec:nablam}

Following~\cite{hhlru}, there is a natural conjecture to make
regarding a combinatorial framework for $\nabla^m(s_\lambda)$ for $m >
1$.  To start, define an \emph{$m$-Dyck path} of length $n$ to be a
lattice path with $n$ north steps and $mn$ east steps that never go
below the line $my = x$.  These reduce to Dyck paths when $m = 1$.  An
\emph{$m$-Dyck sequence}, $g(P) = (g_0,g_1,\ldots,g_{n-1})$, is
defined by setting $g_i$ to be the number of unit squares lying
between the lines $y = i$ and $y = i+1$, to the right of $P$, and
to the left of $my = x$.  An $n$-tuple of nonnegative integers is an $m$-Dyck
sequence for an $m$-Dyck path if and only if $g_0 = 0$ and $g_i \leq
g_{i-1} + m$ for all $1\leq i < n$.

For any given $\lambda$, we define a set $\lndpm{\lambda}$ in a manner
entirely analogous to how we defined $\lndp{\lambda}$.  The primary
difference is that for a pair $(G,R)\in \lndpm{\lambda}$, $G$ is an
\emph{$m$-Dyck configuration}; i.e., each element of $G$ is an
$m$-Dyck sequence.  While $\area(G,R)$ is defined by summing the entries
in $G$ as usual, we generalize the definition of $\dinv(G,R)$ as follows:

\begin{equation}\label{eq:mGRdinv}
\begin{aligned}
  \dinv(G) = \adj(\lambda) &+
  \sumsb{u,v,a,b}\sum_{d=0}^{m-1}\chi(g^{(u)}_a-g^{(v)}_b+d=m) \chi(r^{(u)}_a > r^{(v)}_b)\chi(a\leq b) \\
  &+ \sumsb{u,v,a,b}\sum_{d=0}^{m-1}\chi(1 \leq g^{(u)}_a-g^{(v)}_b+d < m)\chi(a\leq b)\\
  &+ \sumsb{u,v,a,b}\sum_{d=0}^{m-1} \chi(g^{(u)}_a-g^{(v)}_b+d=0)
  \chi(r^{(u)}_a < r^{(v)}_b)\chi(a<b\text{ or }(a=b\text{ and }u<v)).
\end{aligned}
\end{equation}

\begin{conjecture}\label{conj:mMain}
  For any partition $\lambda$ and $m\geq 1$,
  \begin{equation}\label{eq:mfull-conj}
    \nabla^m(s_{\lambda})= \sgn(\lambda)
    \sum_{(G,R)\in \lndpm{\lambda}} t^{\area(G,R)}q^{\dinv(G,R)}x_{R}.
  \end{equation}
\end{conjecture}

We illustrate a typical element of $\lndparb{3}{(2^3)}$ in
Figure~\ref{fig:lndpm}.  As an exercise, the reader can explicitly
write down the elements of $\lndparb{2}{(2,2)}$ labelled with three
$1$'s and a $2$ to compute that
\begin{equation*}
  \langle \nabla^2(s_{2,2}), h_{3,1}\rangle = -(tq)^3(1 + t + t^2 + q + qt + q^2).
\end{equation*}

\myfig{.5}{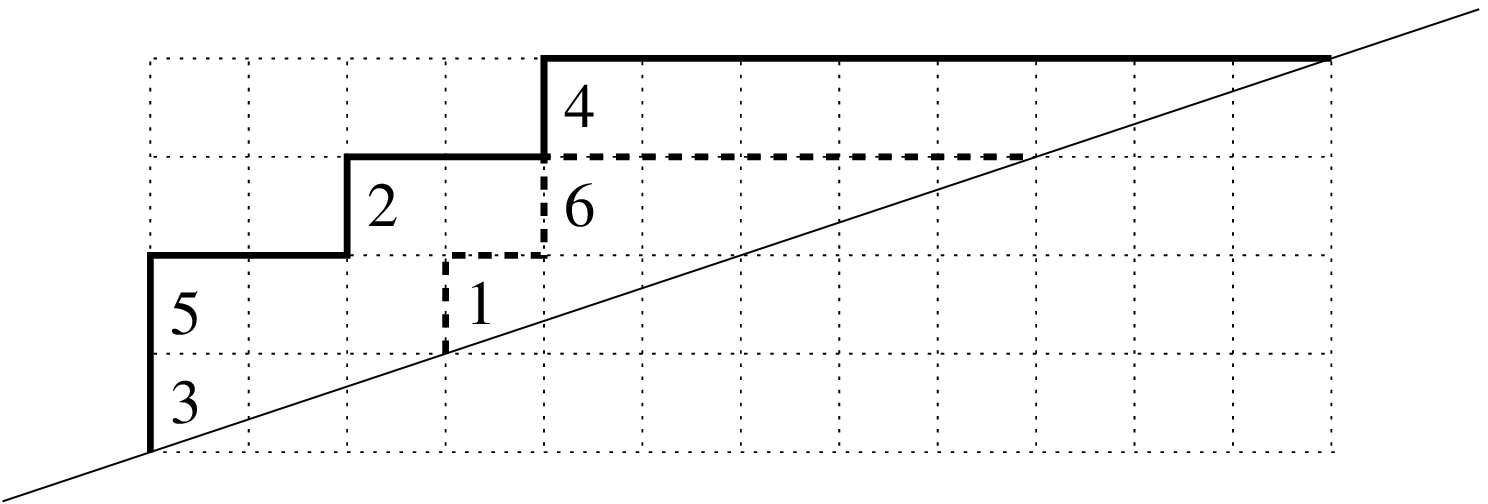}{One term in the coefficient of $s_{1^6}$ in $\nabla^3(s_{2^3})$.}

  We remark that trapezoidal lattice paths and $m$-Dyck paths can be
  treated at the same time by replacing the line $y = x$ with $my = x$ in 
  the definition of trapezoidal paths.  Of course, the $\dinv$ statistic 
  needs to be adjusted accordingly.  See~\cite{thesis,trap} for details.

\section{Proof when $q=1$}
\label{sec:proof}

In this section, we will prove the following specialization
of the main conjectures.
\begin{theorem}
For all partitions $\lambda$,
\begin{eqnarray}
\label{eq:signq}
\scprod{\nabla(s_{\lambda})}{s_{1^n}}|_{q=1} &=&
 \sgn(\lambda)\sum_{G\in NDP_{\lambda}} t^{\area(G)}; \\
\label{eq:fullq}
\nabla(s_{\lambda})|_{q=1} &=&
 \sgn(\lambda)\sum_{(G,R)\in \lndp{\lambda}} t^{\area(G)}x_R.
\end{eqnarray}
\end{theorem}

We remark that Lenart~\cite{lenart} proved a closely related
result that establishes the Schur positivity of $\nabq(s_{\lambda/\nu})$
by expanding the latter polynomials in terms of skew Schur functions.  
Lenart's proof relies heavily on Jacobi-Trudi determinantal formulas.  
We adopt a more combinatorial approach that makes heavy use of 
the inverse Kostka matrix. One benefit of the present method is
that the combinatorial significance of the global sign $\sgn(\lambda)$
is more readily apparent.

\subsection{Specialized Nabla Operator}
\label{subsec:spec-nabla}
The first step is to replace nabla by a more convenient operator.
Let $\Hmuq{\mu}$ denote the image of $\tilde{H}_{\mu}$
under the specialization sending $q$ to $1$; these specialized
Macdonald polynomials form a basis for the $\Q(t)$-vector space
$\Lambda_{\Q(t)}$ of symmetric functions with coefficients in $\Q(t)$.  
Define $\nabq$ to be the unique $\Q(t)$-linear map such that 
$\nabq(\Hmuq{\mu})=t^{n(\mu)}\Hmuq{\mu}$ for all partitions $\mu$. 
It is known that $\nabq$ is a \emph{ring homomorphism} 
(this follows easily from the combinatorial interpretation
of $\tilde{H}_{\mu}$ given in \S\ref{sec:agenda}).
Furthermore, it can be shown that
\[ \nabla(s_{\lambda})|_{q=1}=\nabq(s_{\lambda}). \]
Henceforth we will study the ring homomorphism $\nabq$.

\subsection{Elementary Symmetric Function Expansions}
\label{subsec:e-expansions}
The second step is to study the matrix $\mymat{\nabq}{(e_{\nu})}{(e_{\mu})}$.
Given a Dyck path $\pi$, let $\alpha(\pi)=(\alpha_1(\pi),\alpha_2(\pi),\ldots)$
be the lengths of the vertical columns formed by consecutive north
steps of $\pi$.  Garsia and Haiman proved that
\begin{equation}\label{eq:GHEE}
 \nabq(e_n) = \sum_{\pi\in DP_n} t^{\area(\pi)}e_{\alpha(\pi)} 
\end{equation}
(see Theorem 1.2 in~\cite{remarkable}). Keeping in mind the
combinatorial interpretation of $e_k=s_{1^k}$ in terms
of semistandard tableaux, it is clear that this formula
is equivalent to \eqref{eq:hhlru} when $q=1$.  
Since the ring homomorphism $\nabq$ preserves multiplication, 
we immediately deduce from \eqref{eq:GHEE} that
\begin{equation}\label{eq:nabEE}
 \nabq(e_{\nu})=\sum_{(\pi_1,\pi_2,\ldots)\in DP_{\nu}} 
  t^{\sum \area(\pi_i)}\prod_i e_{\alpha(\pi_i)},
\end{equation}
where $DP_{\nu}$ is the set of all lists of paths $(\pi_1,\pi_2,\ldots)$
such that $\pi_i$ is a Dyck path of order $\nu_i$.
Furthermore, since $\scprod{e_{\xi}}{s_{1^n}}=1$ for all $\xi\vdash n$,
we also deduce that
\begin{equation}\label{eq:nabEsign}
 \scprod{\nabq(e_{\nu})}{s_{1^n}}=\sum_{(\pi_1,\pi_2,\ldots)\in DP_{\nu}} 
  t^{\sum \area(\pi_i)}\qquad(\nu\vdash n).
\end{equation}

\subsection{Transition Matrices}
\label{subsec:trans-matrices}
\begin{figure}
\begin{center}
\epsfig{file=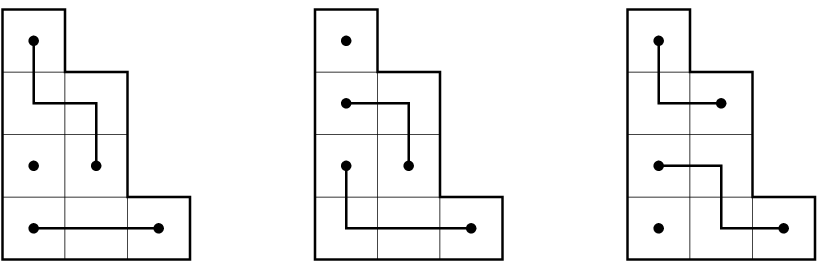}
\caption{Special rim hook tabloids.}
\label{fig:srht}
\end{center}
\end{figure}

The third step is to multiply by suitable transition matrices to change
the input and output bases in the matrix $\mymat{\nabq}{(e_{\mu})}{(e_{\nu})}$.
The relevant transition matrix turns out to be the inverse Kostka matrix,
which we now review. (For more background on transition matrices,
see Section I.6 of~\cite{macbook} or the references~\cite{remtrans,remeg}.) 
Recall that the \emph{Kostka number} $K_{\lambda,\mu}$ is the
number of semistandard tableaux of shape $\lambda$ and content $\mu$.
We have the identities
\begin{eqnarray}
\label{eq:s-to-m}
 s_{\lambda} &=& \sum_{\mu} K_{\lambda,\mu}m_{\mu} \\
\label{eq:h-to-s}
 h_{\mu} &=& \sum_{\lambda} K_{\lambda,\mu}s_{\lambda} \\
\label{eq:e-to-s}
 e_{\mu} &=& \sum_{\lambda} K_{\lambda',\mu}s_{\lambda}
\end{eqnarray}
Letting $K=(K_{\lambda,\mu})$ be the matrix of Kostka numbers, 
these identities assert that $K^t=\mymat{\id}{(s_{\lambda})}{(m_{\mu})}$;
$K=\mymat{\id}{(h_{\mu})}{(s_{\lambda})}$; and
$K=\mymat{\id}{(e_{\mu})}{(s_{\lambda'})}$. Now, let
$K^{-1}=(K_{\lambda,\mu}^{-1})$ be the inverse of the matrix $K$.
The previous identities now read
\begin{eqnarray}
\label{eq:m-to-s}
 m_{\mu} &=& \sum_{\lambda} K_{\mu,\lambda}^{-1}s_{\lambda} \\
\label{eq:s-to-h}
 s_{\lambda} &=& \sum_{\mu} K_{\mu,\lambda}^{-1}h_{\mu} \\
\label{eq:s-to-e}
 s_{\lambda} &=& \sum_{\mu} K_{\mu,\lambda'}^{-1}e_{\mu}.
\end{eqnarray}
Remmel and E\u{g}ecio\u{g}lu discovered the following important
combinatorial interpretation for the entries $K_{\mu,\lambda}^{-1}$
of the inverse Kostka matrix~\cite{remeg}. A \emph{special rim hook
tabloid} of shape $\lambda$ and type $\mu$ is a filling of the
Ferrers diagram $\ferrers{\lambda}$ with rim hooks of length $\mu_i$
that all start in the leftmost column. For example, Figure~\ref{fig:border}
displays one special rim hook tabloid of shape $\lambda'=(5,5,2,1,1)$
and type $(9,5)$. A rim hook spanning $r$ rows has sign $(-1)^{r-1}$.
The sign of a rim hook tabloid is the product of the signs of
all the rim hooks in the tabloid. Let $\srht(\mu,\lambda)$ be
the set of all special rim hook tabloids of shape $\lambda$ and type
$\mu$.  Remmel and E\u{g}ecio\u{g}lu showed that 
$$K_{\mu,\lambda}^{-1}=\sum_{T\in\srht(\mu,\lambda)} \sgn(T).$$ 
For example, we see from Figure~\ref{fig:srht} that
$K_{(4,3,1),(3,2,2,1)}^{-1}=+3$.  
It will be convenient to introduce a modified notion of
the ``type'' of a special rim hook tabloid. Given
a rim hook tabloid $T$, let $\alpha_i(T)$ be the length
of the rim hook that starts $i$ rows from the top
of the diagram for $T$ (for $i=0,1,\ldots$). If
no rim hook starts in row $i$, we let $\alpha_i(T)=0$.
The ordered sequence $\alpha(T)=(\alpha_0(T),\alpha_1(T),\ldots)$
will be called the \emph{total type} of $T$. By dropping
zero entries in $\alpha(T)$ and arranging into decreasing
order, we obtain the type of $T$ (which is a partition).
The total type of the tabloid in Figure~\ref{fig:border}
is $(9,0,0,5,0)$.

\subsection{Intersecting Path Model}
\label{subsec:int-path}

Combining \eqref{eq:s-to-e} with \eqref{eq:nabEE}, we immediately
obtain a combinatorial interpretation for the entries in
the matrix $\mymat{\nabq}{(s_{\lambda})}{(e_{\mu})}$. Using
linearity of $\nabq$, we calculate
\begin{eqnarray*}
\nabq(s_{\lambda}) &=& \sum_{\mu} K_{\mu,\lambda'}^{-1}\nabq(e_{\mu}) \\
&=& \sum_{\mu} K_{\mu,\lambda'}^{-1}\sum_{(\pi_i)\in DP_{\mu}}
  t^{\sum\area(\pi_i)}\prod_i e_{\alpha(\pi_i)}.
\end{eqnarray*}
Here is an explicit combinatorial interpretation of the
right side.  Given $\lambda\vdash n$, we consider all pairs
$(T,\Pi)$ where $T$ is a special rim hook tabloid of shape $\lambda'$ 
and $\Pi=(\pi_i:i\geq 0)$ is a sequence of labelled Dyck paths
such that $\pi_i$ has order $\alpha_i(T)$ for 
$0\leq i<\lambda_1=\ell(\lambda')$.  The sign of such a pair is $\sgn(T)$; 
the $t$-weight of the pair is $\area(\Pi)=\sum_i \area(\pi_i)$;
and the monomial weight is obtained as usual from the labels of $\pi_i$.
The sum of all such signed, weighted objects gives us
$\nabq(s_{\lambda})\in \Q(t)[x_1,x_2,\ldots]$.
Using the fact that $\scprod{e_{\xi}}{s_{1^n}}=1$ for all $\xi$,
we obtain an analogous combinatorial interpretation
for the quantity $\scprod{\nabq(s_{\lambda})}{s_{1^n}}$.
The only difference is that $\Pi$ now consists of unlabelled Dyck paths.

\begin{figure}
\begin{center}
\epsfig{file=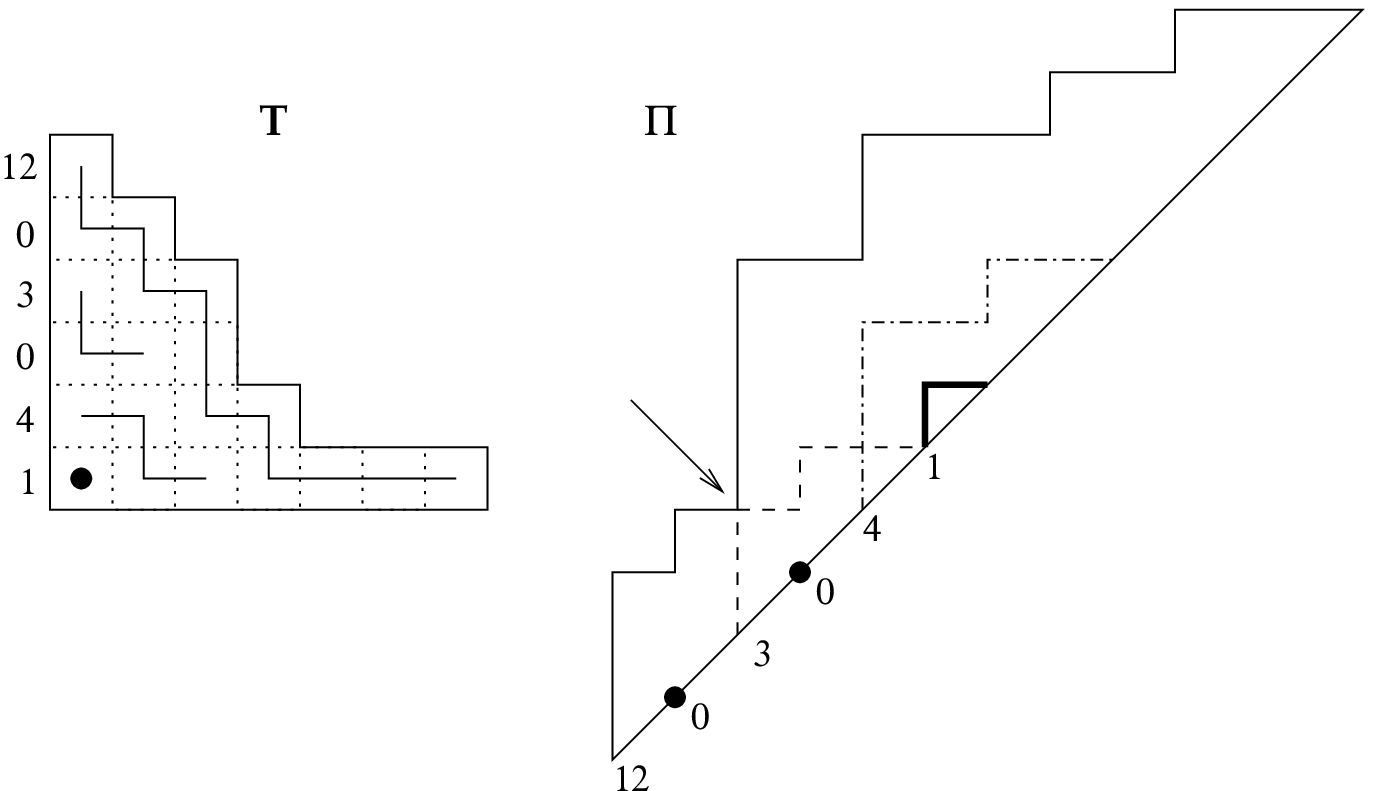}
\caption{Example of a signed object.}
\label{fig:intpath}
\end{center}
\end{figure}

For example, Figure~\ref{fig:intpath} depicts a typical object
contributing to $\scprod{\nabq(s_{(6,5,4,2,1,1,1)})}{s_{1^{20}}}$.  
The sign of this object is $(-1)^7=-1$, and the $t$-weight
is $33+0+2+0+4+0=39$.  As in \S\ref{sec:conjecture}, it is 
convenient to display the paths in $\Pi$ by letting
$\pi_i$ start at $(i,i)$ and end at $(i+\alpha_i(T),i+\alpha_i(T))$.

\subsection{Cancellation for Unlabelled Paths}
\label{subsec:cancel-unlabel}
To complete the proof of the theorem, we define sign-reversing
involutions that cancel objects of opposite sign, and then show that
the fixed points are enumerated by the formulas \eqref{eq:signq}
and \eqref{eq:fullq}. For ease of exposition, we consider the 
unlabelled case first. So far, we have shown that
\begin{equation} 
\label{eq:signq-unc}
 \scprod{\nabq(s_{\lambda})}{s_{1^n}}
  =\sum_{(T,\Pi)} \sgn(T)t^{\area(\Pi)}, 
\end{equation}
where $T$ is any special rim hook tabloid of shape $\lambda'$,
and $\Pi=(\pi_i:0\leq i<\lambda_1)$ is a collection of unlabelled
Dyck paths such that $\pi_i$ has order $\alpha_i(T)$ for all $i$.

\begin{figure}
\begin{center}
\epsfig{file=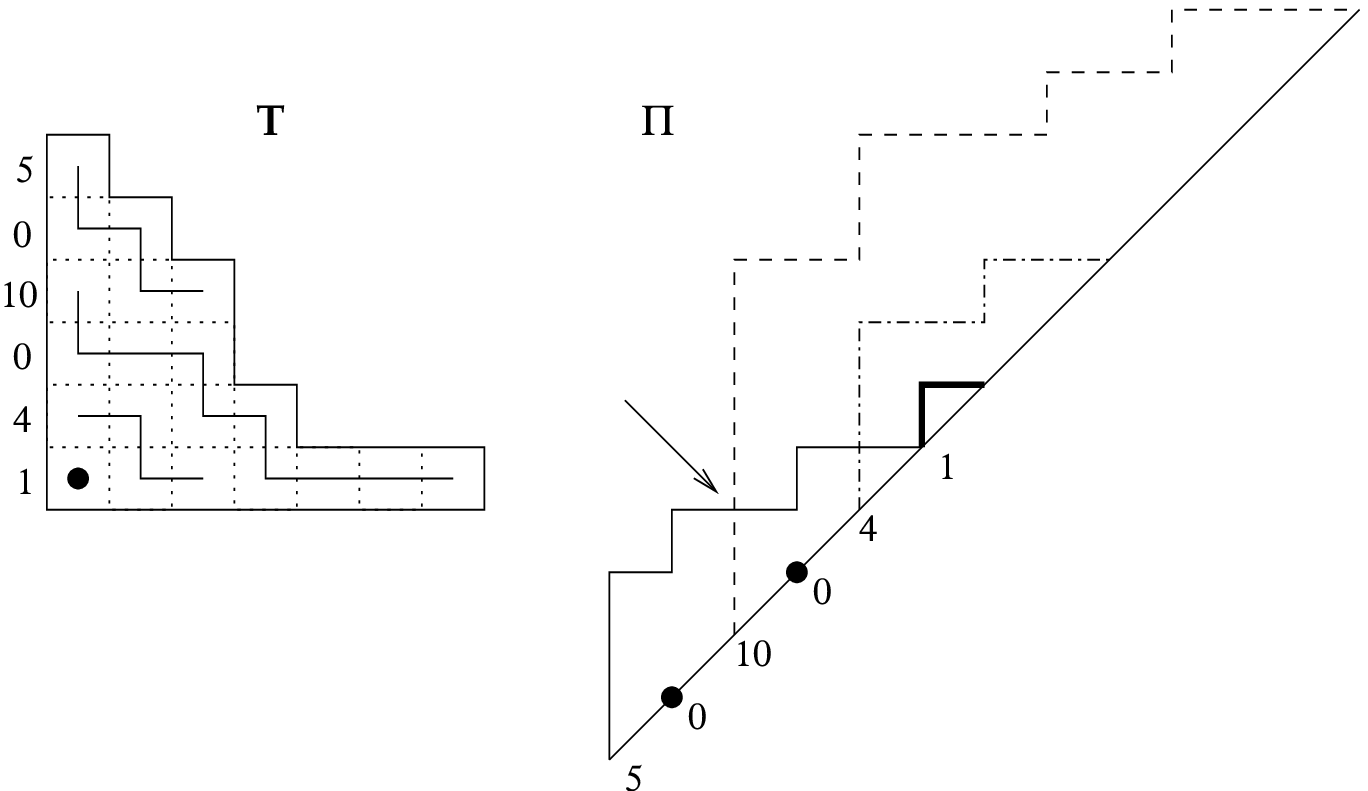}
\caption{Example of tail-switching for paths and tabloids.}
\label{fig:intpath2}
\end{center}
\end{figure}

We can cancel pairs of objects with the same $t$-weight and opposite
signs as follows. Suppose that $(T,\Pi)$ is an object such that
there exist two paths in $\Pi$ that intersect at some vertex.
More precisely, there exist $i<j$ and $(x,y)$ such that
$\pi_i$ and $\pi_j$ both reach $(x,y)$. Among all such
choices of $i,j,x,y$, choose the one such that 
$x$, then $y$, then $i$, then $j$ is minimized.
Write $\pi_i=\beta\gamma$ and $\pi_j=\delta\epsilon$,
where $\beta$ is a path from $(i,i)$ to $(x,y)$,
$\gamma$ is a path from $(x,y)$ to $(i+\alpha_i(T),i+\alpha_i(T))$,
$\delta$ is a path from $(j,j)$ to $(x,y)$,
and $\epsilon$ is a path from $(x,y)$ to $(j+\alpha_j(T),j+\alpha_j(T))$.  
Replace $\pi_i$ by $\beta\epsilon$ and $\pi_j$ by $\delta\gamma$
to get a new list of paths $\Pi'$ with the same earliest
intersection and the same total $t$-weight. Next, consider
the special rim hook tabloid $T$. It is easy to see that
there is a unique way to ``switch the tails'' of the special
rim hooks starting in rows $i$ and $j$ so that the new
special rim hooks in these rows have lengths
$j+\alpha_j(T)-i$ and $i+\alpha_i(T)-j$, respectively.
Furthermore, the sign of the new tabloid $T'$ is opposite
to the sign of $T$. Figure~\ref{fig:intpath2} displays
the object that is matched to the object in Figure~\ref{fig:intpath}
by this process. It is clear from the description that
the map $(T,\Pi)\mapsto (T',\Pi')$ is an involution.  

What do the fixed points of the involution look like?
Clearly, $(T,\Pi)$ is a fixed point iff no two Dyck paths
in $\Pi$ intersect. Because of the way the starting points
of the Dyck paths are arranged along the line $y=x$,
this can only occur if the lengths of the nontrivial
Dyck paths in the list $(\pi_1,\pi_2,\ldots)$ form
a strictly decreasing sequence. In other words, the lengths
of the nonzero rim hooks in $T$ must strictly decrease
reading from top to bottom. One sees easily that this condition
forces $T$ to consist of a succession of ``border hooks''
as described earlier in connection with $\sgn(\lambda)$.
Indeed, we now see that $\sgn(\lambda)$ is simply the
sign of the unique special rim hook tabloid $T$ that
occurs in the objects that are fixed points.
Taking the cancellation and fixed points into account,
we see that the desired result \eqref{eq:signq} follows from
\eqref{eq:signq-unc}.

\subsection{Cancellation for Labelled Paths}
\label{subsec:cancel-label}

We sketch the cancellation for labelled paths while focusing on the
differences with respect to the unlabelled case.  In general, paths
can be rerouted at intersections as in the unlabelled case.  In
particular, any two paths that intersect at the beginning of a common
east step can be cancelled with an object of equal $t$-weight and
opposite sign.  Similarly, if we have an object where one path
intersects the beginning of another (possibly zero-length) path, we
can cancel with an object of equal $t$-weight and opposite sign.  The
label of each step should be envisioned to remain with the individual
step rather than with a particular path.  

We are left to consider the scenario that two paths intersect as in
Figure~\ref{fig:labint}.1 or 2 ($z$ is allowed to be one
in either case).  We claim that each object containing an intersection
of type~1 can be paired with an object containing a region of type~2.  

\myfig{1}{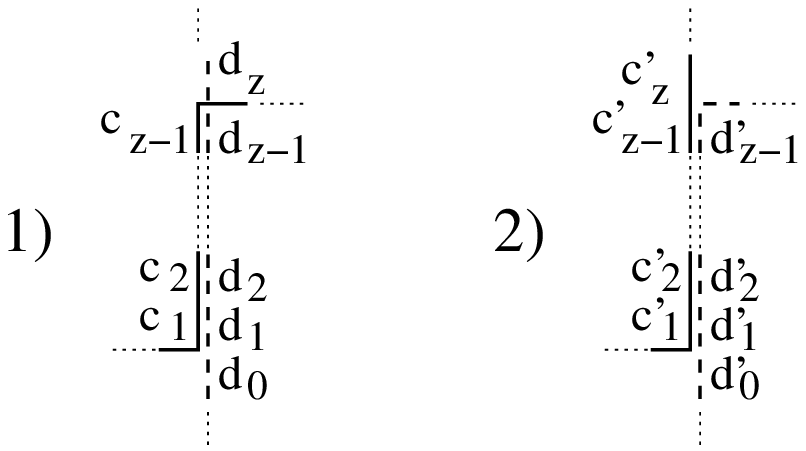}{Cancellation for labelled paths along north steps.}

Assume we have an intersection such as that of Figure~\ref{fig:labint}.1.
Set
\begin{equation*}
  C = \{z\} \cup \{1\leq i< z: c_i > d_{i-1}\} \text{ and }
  D = \{1\} \cup \{1 < i\leq z: d_i > c_{i-1}\}.
\end{equation*}
Suppose that $j\in C\cap D$; choose the smallest such $j$.  Then we
can set $d'_i = d_i$ for $0\leq i< j$, $c'_i = c_i$ for $1\leq i < j$,
$d'_i = c_i$ for $j\leq i < z$ and $c'_i = d_i$ for $j\leq i\leq z$.
This yields a paired object as in Figure~\ref{fig:labint}.2 with the
same $q$- and $t$-weights, but with opposite sign.  This process is
invertible.  So it only remains to show that $C\cap D\neq \emptyset$.

Assume in fact that $C\cap D = \emptyset$.  Hence $1\not\in C$ and
$z\not\in D$.  Let $j$ be as small as possible such that $j\in C$
(such a $j$ must exist since $z\in C$).  It is immediate that $j > 1$
and $j-1\not\in C$.  Therefore, $c_{j-1} \leq d_{j-2} < d_{j-1} <
d_{j}$.  From this we conclude that $j\in D$.  This is a
contradiction.  So, as desired, $C\cap D\neq \emptyset$.

The above argument shows that all configurations with an intersection
such as in Figure~\ref{fig:labint}.1 cancel with an equal-weight
object of the opposite sign.  It remains to characterize the
configurations that do not cancel.  We have just seen that they must
be \emph{weakly nested}.  That is, while any two paths are allowed to
share some north steps, the path that started earlier can never pass
under a path that started later.  In fact, the fixed points are
described by
Conditions~\ref{lndp:cond2},~\ref{lndp:cond3},~\ref{lndp:cond4} along
with
\begin{enumerate}
\item[$4'$.] For all $a$ and $j<k$, if $g_a^{(j)}$ and $g_{a-1}^{(k)}$
  are defined with $g_a^{(j)} = g_{a-1}^{(k)} + 1$, then either
  \begin{enumerate}
  \item $r_a^{(j)} \leq r_{a-1}^{(k)}$, or
  \item $g_a^{(k)}$ is not defined, or
  \item $g_{a-1}^{(j)}$ is not defined, or
  \item $g_a^{(k)}$ and $g_{a-1}^{(j)}$ are both defined with $r_a^{(k)}\leq r_{a-1}^{(j)}$.
  \end{enumerate}
\end{enumerate}

We now show that any fixed point (which must satisfy Condition $4'$),
must actually satisfy Condition~\ref{lndp:cond5} as well.  So suppose
on the contrary we have a fixed point such as in
Figure~\ref{fig:labint}.2 with $a$ and $j<k$ for which $r_a^{(j)} >
r_{a-1}^{(k)}$.  It is easy to check that having $g_a^{(k)}$ or
$g_{a-1}^{(j)}$ not being defined would contradict the fact that we
have a fixed point.  So assume they are, in fact, both defined.  The
only possibility is that $r_a^{(k)}\leq r_{a-1}^{(j)}$.  But then we
give a similar argument to the one above to obtain a contradiction.

\section{Proof Strategy for the Full Conjecture}
\label{sec:agenda}

Most of the known facts about combinatorial interpretations
of the nabla operator (cf. Table~\ref{tab:nabla}) were proved
via long, laborious algebraic manipulations making heavy
use of the machinery of plethystic calculus~\cite{nabla2,nabla3, 
nablaproof1,nablaproof2,schrproof}. In the past, one barrier to 
finding purely combinatorial proofs has been the absence of combinatorial
conjectures that \emph{fully} characterize the action of the
nabla operator. Of course, this barrier is overcome by the
conjecture in \S\ref{sec:conjecture}. Another obstacle to
a combinatorial analysis of nabla was the lack of combinatorial
information about the Macdonald polynomials, which appear in
the definition of nabla. However, recent breakthroughs by
Haglund et al.~\cite{mac-conj,mac-proof,mac-proof2} have resolved this 
difficulty as well.  We review Haglund's combinatorial description of
Macdonald polynomials in the next subsection.  
Combining this material with more inverse Kostka combinatorics,
we then show how to reduce the proof of conjecture~\eqref{eq:full-conj}
to the problem of finding sign-reversing involutions on
certain explicit collections of objects. If this proof
strategy can be completed, it would yield a fully
combinatorial proof of all the conjectured facts about
the nabla operator.

\subsection{Combinatorial Macdonald Polynomials}
\label{subsec:comb-mac}

Let $\mu$ be a fixed integer partition. A \emph{Haglund filling}
of $\mu$ is a function $T:\ferrers{\mu}\rightarrow\N^+$.
Informally, $T$ is a labelling of the cells in the
diagram of $\mu$ with arbitrary positive integers.
The \emph{content monomial of $T$} is
\[ x_T=\prod_{c\in\ferrers{\mu}} x_{T(c)}. \]
Recall that the \emph{major index} of a word
$w=w_1w_2\cdots w_k$ is the sum of the positions
of the descents of $w$, i.e.,  
$\maj(w)=\sum_{i=1}^{k-1} i\chi(w_i>w_{i+1})$.
We define $\maj(T)$ to be the sum of the major
indices of the words obtained by reading the labels
from top to bottom in each column of $T$.

Next, we define the notion of an \emph{inversion triple}.
Consider three cells $c,d,e\in\ferrers{\mu}$ such that
$c$ and $e$ are in the same row (with $c$ to the left of $e$),
and $d$ is the cell immediately below $c$. We call 
$(c,d,e)$ an inversion triple for $T$ iff 
\[ T(c)\leq T(d)<T(e)\mbox{ or }
   T(d)<T(e)<T(c)\mbox{ or }
   T(e)<T(c)\leq T(d). \]
By convention, if $c$ is to the left of $e$ in the bottom
row of $\ferrers{\mu}$, we regard $(c,e)$ as an inversion
triple iff $T(c)>T(e)$. Finally, let $\inv(T)$ be the
total number of inversion triples in $T$.

Haglund~\cite{mac-conj} discovered the formula
\[ \tilde{H}_{\mu}=\sum_{T:\ferrers{\mu}\rightarrow\N^+}
     q^{\inv(T)}t^{\maj(T)}x_T. \]
This formula was proved by Haglund, Haiman, 
and Loehr~\cite{mac-proof,mac-proof2}.
This formula essentially gives the transition matrix
$\mymat{\id}{(\tilde{H}_{\mu})}{(m_{\lambda})}$.
More precisely, writing $x_{\lambda}=\prod_{i\geq 1} x_i^{\lambda_i}$, we have
\begin{equation}\label{eq:H-to-m}
 \tilde{H}_{\mu}=\sum_{\lambda} A(\lambda,\mu)m_{\lambda} 
\mbox{ where }
A(\lambda,\mu)=\sum_{\substack{T:\ferrers{\mu}\rightarrow\N^+ \\
 x_T=x_{\lambda}}} q^{\inv(T)}t^{\maj(T)}.
\end{equation} 

Note that the $\Q(q,t)$-valued coefficient matrix $A=(A_{\lambda,\mu})$
is \emph{invertible}, since the modified Macdonald polynomials
form a basis of $\Lambda$. It follows that the $m_{\lambda}$'s
are the \emph{unique} vectors $v_{\lambda}$ solving the system
of equations $\tilde{H}_{\mu}=\sum_{\lambda} A(\lambda,\mu)v_{\lambda}$.
Now, apply nabla to both sides of \eqref{eq:H-to-m}. We obtain
\begin{equation}\label{eq:nabla-chacz}
 T_{\mu}\tilde{H}_{\mu}=\sum_{\lambda} A(\lambda,\mu)(\nabla(m_{\lambda})).
\end{equation}
Reasoning just as before, it follows that \emph{the vectors 
$\nabla(m_{\lambda})$ constitute the unique solution to the system of equations}
\[ T_{\mu}\tilde{H}_{\mu}=\sum_{\lambda} A(\lambda,\mu)v_{\lambda}. \]
Extracting the coefficient of $m_{\nu}$ on both sides, we see that
\emph{a given indexed family of vectors $(v_{\lambda})$ is equal
to the indexed family $(\nabla(m_{\lambda}))$ iff }
\begin{equation}\label{eq:master1}
 T_{\mu}A(\nu,\mu)=\sum_{\lambda} A(\lambda,\mu)(v_{\lambda}|_{m_{\nu}})
\mbox{ for all partitions }\mu,\nu. 
\end{equation}

\subsection{Combinatorial Formulation of the Problem}
\label{subsec:comb-formulation}

To proceed, we need to have a combinatorial interpretation of the
quantities $\nabla(m_{\lambda})$. One such interpretation follows
immediately from our conjecture for $\nabla(s_{\lambda})$
by using the inverse Kostka matrix again. More precisely, note that
linearity of nabla gives
\begin{eqnarray*}
 \nabla(m_{\lambda}) &=& \sum_{\rho} K_{\lambda,\rho}^{-1}\nabla(s_{\rho})\\
&=& \sum_{\rho}\sum_{T\in\srht(\lambda,\rho)} 
 \sgn(T)\sum_{(G,R)\in\lndp{\rho}}
\sgn(\rho)t^{\area(G,R)}q^{\dinv(G,R)}x_{R},
\end{eqnarray*}
Extracting the coefficient of
$m_{\nu}$, this formula says that the $\nu,\lambda$-entry of
the matrix $\mymat{\nabla}{(m_{\lambda})}{(m_{\nu})}$ is
$$\sum_{\rho}\sum_{T\in\srht(\lambda,\rho)}
\sum_{(G,R)\in\lndp{\rho}:x_R=x_{\nu}}
\sgn(T)\sgn(\rho)t^{\area(G,R)}q^{\dinv(G,R)}.$$

Computer calculations indicate that every entry of
$\mymat{\nabla}{(m_{\nu})}{(m_{\xi})}$ is a polynomial
in $q$ and $t$ with coefficients all of like sign. Yet
the formula just written is a sum of both positive
and negative objects. This indicates that there should
be a sign-reversing, weight-preserving involution on the objects 
just described whose fixed points all have the same sign.
If such an involution could be found, we would have a better
description of the entries of the matrix under consideration.

By the remark at the end of the last subsection, we see that
all the combinatorial formulas for nabla will be proved
if the following ``master identities'' can be verified
for all partitions $\mu,\nu$: 
\begin{equation}\label{eq:master2}
 T_{\mu}A(\nu,\mu)=\sum_{\lambda} A(\lambda,\mu)
\sum_{\rho}\sum_{T\in\srht(\lambda,\rho)}
\sum_{(G,R)\in\lndp(\rho):x_R=x_{\nu}} 
\sgn(T)\sgn(\rho)t^{\area(G,R)}q^{\dinv(G,R)}.
\end{equation} 

Recalling the combinatorial interpretation for the entries of $A$,
we can reformulate the master identities as follows. Fix $\mu,\nu\vdash n$.
On one hand, let $X$ be the set of all tuples $z=(U,\lambda,T,\rho,(G,R))$ 
such that:
\begin{itemize}
\item[(i)] $\lambda$ and $\rho$ are partitions of $n$;
\item[(ii)] $U$ is a Haglund filling of shape $\mu$ and content $\lambda$; 
\item[(iii)] $T$ is a special rim hook tabloid of shape $\rho$
and content $\lambda$;
\item[(iv)] $(G,R)\in \lndp{\rho}$ has content $x_R=x_{\nu}$.
\end{itemize}
The \emph{sign} of $z$ is $\sgn(T)\sgn(\rho)\in\{+1,-1\}$;
the \emph{$t$-weight} of $z$ is $\maj(U)+\area(\Pi)$;
the \emph{$q$-weight} of $z$ is $\inv(U)+\dinv(\Pi)$.

On the other hand, let $X'$ be the set of all Haglund fillings
$U'$ of shape $\mu$ and content $x_{\nu}$. The sign of $U'$
is always positive; the $t$-weight is $n(\mu)+\maj(U')$; the
$q$-weight is $n(\mu')+\inv(U')$. \emph{The master identity \eqref{eq:master2}
holds for $\mu$ and $\nu$ iff there exists a sign-reversing,
weight-preserving involution on $X$ with positive fixed points,
and a weight-preserving bijection between these fixed points and $X'$.}

As a closing remark, we recast the preceding discussion in terms
of matrices. We are essentially trying to prove a matrix identity
of the form $TA=ALB$, where $T$ is the diagonal matrix of scalars
$T_{\mu}$, $A$ is the matrix describing the monomial expansion
of modified Macdonald polynomials, $L$ is the inverse Kostka matrix,
and $B$ is the matrix conjectured to give the monomial
expansion of $\nabla(s_{\lambda})$. Now, $ALB=A(LB)=(AL)B$.
The matrix $LB$ gives the monomial expansion of $\nabla(m_{\lambda})$,
as discussed above, while the matrix $AL$ gives the Schur expansion
of modified Macdonald polynomials. Using special rim hook tabloids,
we can write down collections of \emph{signed}, weighted objects to 
interpret either of the matrix products $LB$ or $AL$.  In each case,
computer evidence (and the known Schur-positivity of Macdonald
polynomials) tells us that we should be able to cancel objects
in these collections to obtain smaller collections of objects
all of like sign. Finding an interpretation for the entries of $AL$
that involves only positive objects is a well-known open problem.
It is likely that this problem (or the analogous problem for $LB$)
will need to be solved first before further progress can be made
in understanding the triple product $ALB$.

\end{document}